\documentclass[11pt,reqno]{amsart} %,fleqn]{amsart}

\usepackage{mathtools}
\mathtoolsset{showonlyrefs}
  
\usepackage{amsmath, amsfonts, amsthm, amssymb,graphicx, color, cite,enumitem}
\usepackage[margin=1.2in]{geometry}

\usepackage{latexsym,hyperref}

\usepackage{hyperref}
\hypersetup{hidelinks}

\newtheorem{Theorem}{Theorem}[section]

\theoremstyle{definition}

\theoremstyle{remark}
\newtheorem{Remark}{Remark}[section]

\numberwithin{equation}{section}
\allowdisplaybreaks

\renewcommand{\u}{{\bf u}}

\def\f{\frac}

\def\hf1{^\f{1}{1-\xi^2}}

\def\be{\begin{equation}}
\def\en{\end{equation}}
\def\bs{\begin{split}}
\def\es{\end{split}}
\def\ba{\begin{align}}
\def\ea{\end{align}}

\author[A. L. Mazzucato]{Anna L. Mazzucato}
\address{Department of Mathematics,  Penn State University, University Park, PA 16802,
USA.}
\email{alm24@psu.edu}

\author[D. Wang]{Dehua Wang}
\address{Department of Mathematics, University of Pittsburgh, Pittsburgh, PA 15260, USA.}
\email{dhwang@pitt.edu}

\author[W. Wei]{Wei Wei}
\address{Department of Mathematics, University of Pittsburgh, Pittsburgh, PA 15260, USA.}
\email{wew144@pitt.edu}

\title[Vanishing viscosity with inflow/outflow]
{On the vanishing viscosity limit for incompressible flows with inflow/outflow boundary conditions}

\date{\today}

\begin{document}
\begin{abstract}
We study the vanishing viscosity limit for the incompressible Navier-Stokes equations (NSE) in a general bounded domain with inflow-outflow boundary conditions. Extending the work of Gie, Hamouda, and Temam ({\em Netw. Heterog. Media} {\bf 7}, 2012) and also of Lombardo and Sammartino ({\em SIAM J. Math. Anal.} {\bf 33}, 2001), we allow for a general injection and suction angle, as long as it is bounded away from zero. We rigorously establish the convergence of NSE solutions to those of the Euler equations (EE) as viscosity vanishes in the energy norm. We prove interior convergence in both the $L^2$ and the Sobolev $H^1$ norms at the same rates as in the case of injection/suction normal to the boundary. The proof relies on the construction of boundary layer correctors via Prandtl-type equations and a higher-order asymptotic expansion that improves the convergence rate. 

\end{abstract}

\maketitle

% section 1. introduction
\section{Introduction}

This work concerns the behavior of incompressible viscous fluids as viscosity vanishes in bounded domains in two and three space dimensions, in the presence of permeable boundaries.

Viscous incompressible flows are modeled by the {\em Navier-Stokes equations} (NSE), while the {\em Euler equations} (EE) describe the motion of incompressible inviscid fluids. Formally, the Navier-Stokes equations become the Euler equations as the viscosity parameter vanishes, but this limit is a singular limit that is generally poorly understood analytically in the presence of boundaries. Experimentally, it is well established that in a thin region adjacent to the boundary - the so-called {\em viscous boundary layer} - viscosity cannot be neglected no matter how small it is. Indeed, friction at walls generates vorticity that can generate instabilities and, ultimately, onset of turbulence.

Mathematically, there are two related, but distinct problems. One is the convergence of Navier-Stokes solutions to Euler solutions as viscosity vanishes. We refer to this limit as the {\em vanishing-viscosity} or {\em zero-viscosity limit}. Mathematically, this limit holds if a sequence of Navier-Stokes solutions converges to a sequence of Euler solutions strongly in the energy norm, that is, uniformly in time in the Lebesgue space $L^2$. 
Because of the ill-posedness of the Euler equations in the context of weak solutions, in the zero-viscosity limit, one typically assumes that the limit Euler solution is a sufficiently regular solution, while the Navier-Stokes solutions are typically Leray-Hopf weak solutions, so that the energy inequality holds.
The other problem is to model the behavior of slightly viscous fluids in the boundary layer. It goes back to the work of L. Prandtl to derive effective equations via asymptotic analysis based on scaling with respect to the viscosity parameter, now called {\em Prandtl equations}. These are degenerate equations that are generally more challenging to study than the NSE and EE, since they can blow up in finite time or exhibit ill-posedness. More refined models have also been proposed, such as the {\em triple deck} and {\em interactive boundary layer} models. (We refer the reader to the monographs \cite{Gie2018,SchlichtingBook2017,QinEtAl2024,CousteixMauss2007}, the survey paper \cite{MaekawaMazzucato2018}, and references therein.)

Formally, the zero-viscosity limit is a singular perturbation problem, since the viscosity coefficient characterizes the highest-derivative in the equation. Even so, in the absence of physical boundaries, when the fluid occupies the entire space, or under periodic boundary conditions, the vanishing viscosity limit is relatively well understood (see e.g. \cite{Constantin1995,Kato1972,Swann1971,Yudovitch1963,Masmoudi2007}, and also \cite{BardosEtAl2012} in the context of irregular Euler solutions). It is also possible to establish the vanishing viscosity limit in situations where the domain of the viscous fluids changes in such a way to become negligible in the limit, as in the case of vanishing obstacles \cite{IftimieEtAl2009,IftimieKelliher2009,LacaveMazzucato2016,Viana2021}, or expanding domains \cite{KelliherEtAl2009}.

By contrast, the presence of physical boundaries that exercise friction on the movement of the viscous fluid introduces significant challenges. Analytically, the mismatch between the boundary conditions for Navier-Stokes and Euler equations can allow for rapid changes of the velocity in the boundary layer which can lead to growth of derivatives.  For impermeable walls, under the classical {\em no-slip} boundary conditions for Navier-Stokes equations,  the viscous velocity vanishes on the boundary, and strong boundary layers form due to an incompatibility with Euler solutions, which satisfies only the condition of no normal flow, the so-called {\em no-penetration} condition. This setting leads to a strongly singular perturbation and remains largely open. In the unsteady case, 
 the limit is known to hold in three cases essentially: when the problem is suitably linearized \cite{Temam1997,TemamWang1995,TemamWang1996.a,TemamWang1996.b,TemamWang1997.b,Hamouda2008,LombardoEtAl2001}, when the data is either analytic or monotonic \cite{Oleinik1999,NguyenNguyen2024,KukavicaEtAl2024,CaflischSammartino1997,SammartinoCaflisch1996,SammartinoCaflisch1998.a,SammartinoCaflisch1998.b,Sammartino1997,Maekawa2014,NguyenSquare2024,NguyenSquare2018,BardosEtAl2022,Wang2020,WangEtAl2020,KukavicaEtAl2014} (but see \cite{FeiEtAl2018} for convergence with non-analytic data and a control on vorticity at the boundary) or under special symmetry on the flow \cite{Matsui1994,BonaWu2002,MazzucatoEtAl2011,GieEtAl2019,MazzucatoTaylor2008,MazzucatoTaylor2011,LopesEtAl2008.a,LopesEtAl2008.b,HanEtAl2012}.
There are also several conditional convergence results, starting with the seminal work by Kato \cite{Kato1984}. We mention in particular the works \cite{Wang2001,Kelliher2007,Kelliher2017,Kelliher2008,Kelliher2023,ChenEtAl2022,BardosEtAl2019,BardosTiti2007,BardosTiti2013,ConstantinEtAl2017,ConstantinEtAl2018,ConstantinEtAl2019,Shuai2025,Sueur2012}. An interesting unconditional result concerns the maximal separation rate in the vanishing viscosity limit \cite{VasseurYang2023,VasseurYang2024}, while bounds on the viscous energy dissipation rate with constant injection and suction were obtained in \cite{DoeringEtAl2000}.

Another physically important boundary condition for viscous flows in impermeable domains is the {\em Navier slip-with-friction} condition, which allows the flow to slip at the boundary in a manner proportional to the tangential stress. Under these conditions, it is possible to show the validity of the vanishing viscosity solutions under some assumptions on the friction coefficient and on the boundary. Since these conditions are not the focus of this work, we only mention some references among the vast literature \cite{ClopeauEtAl1998,LopesPlanas2005,Rusin2006,MasmoudiRousset2012,WangEtAl2012,XiaoXin2007,Zhong2017,Liu2024,Mucha2005,NeustupaPenel2011,BerselliSpirito2012,GerardVaretEtAl2018,DrivasNguyen2019,TaoEtAl2020,Nguyen2020,daVeigaCrispo2011,daVeigaCrispo2011,GieKelliherNavier2012,Gie2013} (see also the survey \cite{daVeigaCrispo2023} and references therein).

This paper focuses on the vanishing viscosity limit and viscous layers when rigid walls are {\em permeable}, allowing for non-zero normal flow for both viscous and inviscid fluids, specifically under {\em injection} and {\em suction} at the boundary. This setting is especially important in applications, such as in modeling flows in sections of pipes and channels, flow through (stiff) biological membranes, or in moving domains, for instance. Experimentally, but also computationally and via asymptotic analysis, it is observed that a sufficiently strong injection and suction rate can stabilize the viscous boundary layer (see e.g. \cite{IshakEtAl2007,Magyari2000}. Therefore, it is expected that the vanishing viscosity limit holds with these boundary conditions that make the boundary non-characteristic for the problem. However, while injection and suction problem has been well studied in the context of compressible flows, such as in jets or flows above obstacles, it has received less attention in the case of incompressible flows. In fact, there are analytic challenges in dealing with injection and suction for inviscid fluids. For the Navier-Stokes equations, there is no conceptual difficulties in imposing non-homogeneous and non-tangential velocities at the boundary. However, for the Euler equations, specifying only a (non-zero) normal component of the velocity, with the appropriate sign for inflow/outflow, does not seem sufficient generically. On the one hand, analytically, it leads to a potential derivative loss at the boundary, at least in Sobolev spaces. On the other hand, vorticity is produced by normal flow at the boundary and is carried inside the fluid domain at the inflow. Therefore, the vorticity or pressure must be determined at the inflow. By restricting the equations to the inflow boundary, a formally determined problem can be obtained by also imposing the {\em full} velocity at the inflow (we refer the reader to a discussion of this delicate point in \cite{Mazzucato2025}). Imposing the full inflow velocity or the normal velocity and the pressure change at the inflow are the more physical conditions.

One of the first results concerning the Euler equations with inflow/outflow is due to Yudovich \cite{Yudovich1962,Yudovich1964}, who proved global existence and uniqueness of strong solutions, assuming the vorticity is given at  inflow and it is generated by a compatible vector field (we mention also the recent works \cite{Bravin2022,BravinSueur2022,NoisetteSueur2023} dealing with sources and sinks in this setting, and the article \cite{ChemetovAntotsev2008}, where a friction boundary conditions is imposed at the inflow). Bardos and Ghidaglia constructed weak solutions when the head pressure, i.e., the Bernoulli function, is specified at inflow along with the normal component of the velocity \cite{BardosGhidaglia1998}. 
The well-posedness of the Euler equations with compatibility condition only on the data at the initial time can be found in \cite{Antontsev1990}, where the short-time existence of strong solutions to the Euler equations in the H\"older space $C^{1,\alpha}$ was obtained specifying the velocity at inflow (see \cite{Petcu2006} for further well-posedness results). This result was recently generalized by one of the authors of this article and her collaborators to well-posedness in the H\"older space $C^{N,\alpha}$ for any $N\in \mathbb{N}$ \cite{Mazzucato2025,GieEtAl2023}. Compatibility conditions are not needed with analytic data, since well-posedness follows by imposing only the normal component of the velocity at the boundary, and long-time existence can be obtained in the two-dimensional case with data chosen suitably small \cite{KukavicaEtAl2025}. 

Contrary to the inviscid case, the theory of both weak and strong solutions for the Navier-Stokes equations is well developed, assuming inhomogeneous {\em slip}, i.e., Dirichlet, boundary conditions. Mirroring the no-slip case, local-in-time well-posedness of strong solutions holds in various function spaces (see \cite{FursikovEtAl2008,FarwigEtAl2006,FarwigEtAl2017,FarwigEtAl2010,Raymond2007,ChangJin2020} and references therein). Furthermore, in the viscous case it is possible to treat the case when the outflow, inflow, and impermeable parts of the boundary touch, at least as long as the boundary data is continuous. However, corner-type singularities are still possible and have been considered \cite{Zajaczkowski1987.a,Zajaczkowski1987.b}. This setting is especially important in applications,  for example in modeling flows in sections of pipes and channels (see also \cite{Collini1986}). For the Euler equations instead, we have to assume that the inflow/outflow and impermeable boundary components are closed and disjoint; otherwise there is a potential jump in the boundary velocity,  which we only fully control at inflow, where these components meet.

Assuming the existence of a sufficiently regular Euler solution with injection and suction conditions at the boundary, it is possible to rigorously study the zero-viscosity limit and establish a Prandtl-type asymptotic expansion for the viscous solution in a boundary layer.  This result rigorously justifies the observed phenomenon that injection and suction stabilize the layer. In fact, at inflow the Euler and Navier-Stokes equations have compatible data so that correctors are not needed, at least at first order, to obtain convergence, while at outflow the Prandtl boundary layer is of size proportional to the viscosity and not its square root as in the no-slip case, hence much weaker.

Alekseenko \cite{Alekseenko1994} first studied $L^2$ convergence of Navier-Stokes solutions to Euler solutions in domains characterized by an inlet, an outlet, and lateral surfaces forming curvilinear dihedral angles. However, this work does not address the boundary layer behavior or stronger convergence in the interior. Later, Temam and Wang \cite{Temam2002} proved interior convergence in the Sobolev space $H^1$ in a three-dimensional periodic channel with injection and suction that is normal to the boundary, while  in two dimensions they also established uniform convergence. Gie, Hamouda, and Temam \cite{Gie2012} extended these results to curved domains in $\mathbb{R}^3$ again with inflow-outflow normal to the boundary, and improved the convergence rate by using higher-order asymptotic expansions.  Lombardo and Sammartino \cite{Lombardo2001} considered the linearization of the Navier-Stokes equations around a steady shear flow, an Oseen-type equation, in a two-dimensional periodic channel, but allowed for injection and suction at an oblique angle to the walls. They established the validity of the zero-viscosity limit and proved $H^1$ convergence in the interior. (See also related results in \cite{Gulas2025, Temam2000, chemetovCipriano2014}.)

In this work we generalize the results of \cite{Gie2012,Lombardo2001}. We consider the vanishing viscosity problem in general bounded domains enclosed by smooth inlet and outlet surfaces. Unlike previous studies, we allow oblique inflow-outflow directions that are not necessarily aligned with the boundary normal. Within this more general geometric framework, we establish the vanishing viscosity limit and recover the convergence rates of \cite{Gie2012}
(see Theorem \ref{t:main}). We employ suitable curvilinear coordinates in the boundary layer and construct layer correctors to control both the mismatch in the tangential component of the velocity at outflow, as well as boundary derivatives and the pressure. Although already mentioned in \cite{Gie2012},  we explicitly allow for the presence of umbilical point on the boundary surface in three space dimensions, which prevents us from using a global coordinate system based on the boundary normal and principal curvatures on the surface. Instead, we utilize local boundary normal coordinates and a partition of unity. 

Our main contribution is a unified approach to the vanishing viscosity limit with inflow/outflow, which leads to establishing a very general result, including higher-order convergence in the interior and convergence in $H^1$ of the corrected Euler solution for sufficiently regular and compatible data. In turn, these results give indirectly bounds for the vorticity production at the boundary (vorticity may be produced even at zero viscosity by the outflow/inflow). Our results, in particular, only require that injection and suction do not degenerate along the boundary. In comparison, Yang, Martinez, Vasseur, and one of the authors   showed in \cite{YangetAlArXiV2024} that, without compatibility conditions, the rate of possible separation between the Navier-Stokes solution and the Euler solution  is consistent with that in the impermeable case and anomalous energy dissipation holds, if suction is sufficiently strong. %\cite{YangetAlArXiV2024}.

The paper is organized as follows. Section \ref{s:prelim} introduces the problem, states the main result, and develops the necessary differential operators in curvilinear coordinates. Section \ref{s:Prandtl}  contains a derivation of the Prandtl equation in our setting and the construction of the leading-order boundary layer corrector. Section \ref{s:Sobolev} establishes convergence in the $L^2$ and $H^1$ norms. Finally, Section \ref{s:correctors} refines the asymptotic analysis to give higher-order boundary layer correctors, improving the convergence rate.  

Throughout, we follow the standard notation. For instance, $L^p$, $1\leq p \leq \infty$, denotes the space of $p$th-integrable functions, while $H^s$, $s\geq 0$, is the $L^2$-based Sobolev space of order $s$. We will use the shorthand notation $L^p\,X$ to denote the space of functions $f:[0,T)\to X$ that are $p$th integrable in time with values in the Banach space $X$. Furthermore, $C$ indicates a generic constant that may change from line to line. We follow the notation typical of asymptotic analysis and denote the viscosity coefficient, our small parameter, with $\epsilon$ instead of the more common $\nu$ in fluid mechanics. In particular, with $f\approx g$ we mean that $f=g$ up to terms that are formally $o(1)$ in $\epsilon$.

% section 2. formulation and notation
\section{Formulation of the Problem and Preliminaries} \label{s:prelim}

In this section, we formalize the vanishing viscosity limit in the context of the incompressible Navier-Stokes equations  in a bounded domain with inflow/outflow and introduce some needed geometric concepts. For notational ease, we refer to the Navier-Stokes equations also as NSE, while we refer to the Euler equations also as EE. We confine ourselves to the three-dimensional case, but similar results hold in two space dimensions. For the vanishing viscosity limit, in general the same difficulties exist in two or three dimensions, although the geometry is simplified in the two-dimensional case.

In Subsection \ref{setting}, we discuss the NSE and a domain in $\mathbb{R}^3$, bounded by two disjoint surfaces representing the inflow and outflow boundaries.  The inflow-outflow directions are allowed to be oblique, i.e., not necessarily aligned with the boundary normal. In Subsection \ref{coord}, we introduce a curvilinear coordinate system adapted to the geometry near the outflow boundary.

% subsection 2.1. settings for the vanishing viscosity problem
\subsection{The vanishing viscosity limit with inflow/outflow} \label{setting}

Let $\Omega\subset \mathbb{R}^3$ be a regular domain with boundary $\Gamma = \Gamma_- \cup \Gamma_+$, where $\Gamma_-$ and $\Gamma_+$ are two disjoint closed surfaces representing the outflow and inflow portions of the boundary, respectively.  

We consider the NSE with non-homogeneous boundary data in the space-time cylinder $\Omega\times [0,T)$:
\begin{equation}
\begin{cases}
\partial_t u^\epsilon - \epsilon \Delta u^\epsilon + (u^\epsilon \cdot \nabla) u^\epsilon + \nabla p^\epsilon = f,\\
\operatorname{div} u^\epsilon = 0,\\
\gamma_- u^\epsilon = \alpha,\\
\gamma_+ u^\epsilon = \beta,\\
u^\epsilon|_{t=0} = u_0,
\end{cases}
\end{equation}
where $\u^\epsilon$ represents the velocity of the fluid, $p^\epsilon$ is the pressure, $\epsilon$ is the viscosity coefficient, $f$ are given external forces, $u_0$ is the initial condition, $\alpha$, $\beta$ represent the prescribed velocity on the two components of the boundary with $\gamma_\pm$ the trace operator on $\Gamma_\pm$.
The boundary data are chosen to satisfy the inflow/outflow condition, that is:
\begin{equation}\label{BC compatibility 1}
\alpha \cdot n > 0, \quad \beta \cdot n < 0,
\end{equation}
where $n$ denotes the outward unit normal on $\Gamma$. 
The incompressibility condition gives the additional constraint on the boundary data: 
\begin{equation}\label{BC compatibility 2}
\int_{\Gamma_-} \alpha \cdot n \, dS + \int_{\Gamma_+} \beta \cdot n \, dS = 0,
\end{equation}
with $dS$ infinitesimal surface area.
We assume that the data is independent of the viscosity, but this condition can be relaxed. For example, the initial condition could be taken instead of the form $u^\epsilon(0)=u_0^\epsilon$ with $u^\epsilon_0$ converging to $u_0$ as $\epsilon \to 0$ in a suitable way.

We may assume that the boundary data are induced by a background flow $\mathcal{U}$ (see e.g. \cite{Temam1979}), which solves the following Stokes' problem:
\begin{equation}
\begin{cases}
- \Delta \mathcal{U} + \nabla p = 0,\\
\operatorname{div} \mathcal{U} = 0,\\
\gamma_- \mathcal{U} = \alpha,\\
\gamma_+ \mathcal{U} = \beta,
\end{cases}
\end{equation}
for a certain pressure. 

By setting $\epsilon = 0$, the NSE reduces formally to the EE:
\begin{equation}
\begin{cases}
\partial_t u^0 + u^0 \cdot \nabla u^0 + \nabla p^0 = f,\\
\operatorname{div} u^0 = 0,\\
\gamma_- u^0 \cdot n = \alpha \cdot n,\\
\gamma_+ u^0 = \beta,\\
u^0|_{t=0} = u_0,
\end{cases}
\end{equation}
where only the normal component of the velocity is specified on $\Gamma_-$. The mismatch in the tangential components between the NSE and EE gives rise to a potentially strong boundary layer near $\Gamma-$ in the vanishing viscosity limit.

We assume all data are  smooth:
\begin{equation}\label{regularity}
\begin{split}
&f \in C^\infty([0,T] \times \Omega), u_0 \in C^\infty(\Omega),\\
&\Gamma \in C^\infty, \alpha \in C^\infty([0,T] \times \Gamma_-), \beta \in C^\infty([0,T] \times \Gamma_+),
\end{split}
\end{equation}
and  satisfy the compatibility conditions stated in \cite{Mazzucato2025}, ensuring that both the EE and NSE possess unique smooth solutions for a short time. 
We will not concern ourselves with the optimal regularity for the data, which is not the focus in the vanishing viscosity limit and to streamline presentation. We observe that, by standard elliptic estimates and linearity,   $\mathcal{U}, p \in C^\infty([0,T] \times \Omega)$ provided that $\alpha \in C^\infty([0,T] \times \Gamma_-)$ and $\beta \in C^\infty([0,T] \times \Gamma_+)$.

Since we employ correctors, we require the additional compatibility condition:
\begin{equation}\label{compatibility0}
   \gamma_- u^0 |_{t=0} \cdot \tau = \gamma_-u_0 \cdot \tau = \alpha |_{t=0} \cdot \tau.
\end{equation}
Furthermore, to establish interior convergence in the vanishing viscosity limit, we impose that
\begin{equation}\label{compatibility}
%\begin{cases}
\alpha \cdot n \geq \alpha_0, \text{ for some } \alpha_0 > 0.
%\end{cases}
\end{equation}

It is a standard approach to shift the solution so that it satisfies homogeneous boundary conditions. To this end, we define the perturbations   $v^\epsilon : = u^\epsilon-\mathcal{U}$ and $v^0 : =  u^0-\mathcal{U}$, which then solve the homogenized systems:
\begin{equation}\label{NSE}
\begin{cases}
\partial_t v^\epsilon - \epsilon \Delta v^\epsilon + \mathcal{U} \cdot \nabla v^\epsilon + v^\epsilon \cdot \nabla \mathcal{U} + v^\epsilon \cdot \nabla v^\epsilon + \nabla p^\epsilon = f + \epsilon \Delta \mathcal{U},\\
\operatorname{div} v^\epsilon = 0,\\
\gamma v^\epsilon = 0,\\
v^\epsilon|_{t=0} = v_0,
\end{cases}
\end{equation}
and, respectively:
\begin{equation}\label{EE}
\begin{cases}
\partial_t v^0 + \mathcal{U} \cdot \nabla v^0 + v^0 \cdot \nabla \mathcal{U} + v^0 \cdot \nabla v^0 + \nabla p^0 = f,\\
\operatorname{div} v^0 = 0,\\
\gamma_- v^0 \cdot n = 0,\\
\gamma_+ v^0 = 0,\\
v^0|_{t=0} = v_0,
\end{cases}
\end{equation}
where the term $- (\partial_t \mathcal{U} + \mathcal{U} \cdot \nabla \mathcal{U})$ has been absorbed into $f$. From the compatibility condition \eqref{compatibility0}, we deduce:
\begin{equation}\label{compatibility 2}
\gamma_- v^0|_{t=0} \cdot \tau = \gamma_- v_0 \cdot \tau = 0.
\end{equation}

To achieve higher-order interior convergence, we acquire an additional compatibility condition
\begin{equation}\label{compatibility 3}
\gamma_- (\partial_t u^0)|_{t=0} \cdot \tau = (\partial_t \alpha)|_{t=0} \cdot \tau,
\end{equation}
which is equivalently expressed in terms of $v^0$ as
\begin{equation}\label{compatibility 4}
\gamma_- (\partial_t v^0)|_{t=0} \cdot \tau = 0.
\end{equation}

The main result of this paper is the following theorem.

%\AnnaComment{I don't think that we can include the estimate on the time $T$, because the Euler solution may exists for a shorter time. I have included a remark}
\begin{Theorem} \label{t:main}
Assume the regularity and compatibility conditions \eqref{regularity} and \eqref{compatibility 2}. Then there exists a sufficiently small $T>0$ and a constant $C>0$ independent of $\epsilon$ such that
\begin{equation}
\begin{cases}
\|v^\epsilon - (v^0 + \varphi^0)\|_{L^\infty L^2} \leq C \epsilon,\\
\|v^\epsilon - (v^0 + \varphi^0)\|_{L^2 H^1} \leq C \epsilon^{1/2},\\
\end{cases}
\end{equation}
and
\begin{equation}
\begin{cases}
\|v^\epsilon - (v^0 + \varphi^0) - \epsilon (v^1 + \varphi^1)\|_{L^\infty L^2} \leq C \epsilon^2,\\
\|v^\epsilon - (v^0 + \varphi^0) - \epsilon (v^1 + \varphi^1)\|_{L^2 H^1} \leq C \epsilon^{3/2},\\
\end{cases}
\end{equation}
provided the additional compatibility condition \eqref{compatibility 4} holds.
The correctors $\varphi^0$, $v^1$, $\varphi^1$ are defined in \eqref{BL formula}, \eqref{bulk 1}, \eqref{higher order potential}, and \eqref{higher order BL formula}.
\end{Theorem}

\begin{Remark}
 If the Euler solution persists in time, such as in the case of two space dimensions and sufficiently small analytic data, then a lower bound of the time $T$ is given in terms of the estimates in \eqref{small time}.
\end{Remark}

% subsection 2.2. differential operators in curvilinear coordinates
\subsection{Differential operators in curvilinear coordinates}\label{coord}

We adopt the convention that subscripts and superscripts indicate covariant and contravariant components or basis vectors, respectively.

Let $\mathcal{A}(\Gamma_-) = \{\Phi_{(k)}: I_2 \to \Gamma_{(k)}\}_{i=1}^N$ 
be an atlas of $\Gamma_-$, where $I_2 = (0,1)^2 \subset \mathbb{R}^2$, the collection 
$\{\Gamma_{(k)}\}_{i=1}^N$ forms an open cover of $\Gamma_-$, 
and each $\Phi_{(k)}$ is a $C^\infty$ local chart.

We employ {\em boundary normal} coordinates. For the reader's sake, we briefly recall how these are constructed. 
Let $\delta \in (0, \operatorname{dist}(\Gamma_-, \Gamma_+))$ be small enough.  For each $k$, define a local diffeomorphism $\bar{\Phi}_{(k)}: I \rightarrow H_{(k)}$ by
\begin{equation}
\bar{\Phi}_{(k)}(\xi) = \Phi_{(k)}(\xi') - \xi_3 n,
\end{equation}
where $I = I_2 \times [0,\delta]$, $\xi=(\xi',\xi_3)$ denotes a point in $I$, and $H_{(k)} = \bar{\Phi}_{(k)}(I)\subset \overline{\Omega}$. In particular $\xi_3$ gives the distance of the point $x\in H_k$ to $\partial \Omega$. This
identification induces a curvilinear coordinate system in $H_{(k)}$ with covariant basis vectors $e_{(k),i} = \frac{\partial}{\partial \xi_i} {\Phi}_{(k)} - \xi_3 \frac{\partial n}{\partial \xi_i}$, for $i=1,2$, and $e_{(k),3} = - n$. Since $e_{(k),i} \cdot e_{(k),3} = 0$, for $i=1,2$, and $\|g_{(k),3}\| = 1$, the metric tensor takes the form:
\begin{equation}
(E_{(k),ij}) = 
\begin{pmatrix}
E_{(k),11} & E_{(k),12} & 0 \\
E_{(k),21} & E_{(k),22} & 0 \\
0 & 0 & 1
\end{pmatrix},
\end{equation}
with inverse
\begin{equation}
(E^{ij}_{(k)}) = (E_{(k),ij})^{-1} =
\begin{pmatrix}
E^{11}_{(k)} & E^{12}_{(k)} & 0 \\
E^{21}_{(k)} & E^{22}_{(k)} & 0 \\
0 & 0 & 1
\end{pmatrix}.
\end{equation}
We denote $V_{(k)} = |(E_{(k),ij})|^2$.
%\AnnaComment{I would prefer to use $e_{(k)}$ for the basis vectors.}

We now summarize several differential operators expressed in this coordinate system (the subscript $(k)$ is omitted for clarity). Let f be a smooth scalar function, and F, G smooth vector fields. Then:\\
\begin{enumerate}
 \item Gradient of a scalar field:
\begin{equation}\label{gradient}
%\begin{split}
\nabla f = \frac{\partial f}{\partial \xi_j} E^{ij} e_i
= \frac{\partial f}{\partial \xi_3} e_3 + \sum_{i,j=1,2} \frac{\partial f}{\partial \xi_j} E^{ij} e_i;
%\end{split}
\end{equation}
\item  Directional derivative:
\begin{equation}\label{directional derivative}
\begin{split}
\nabla_F G &= F^j \left(\frac{\partial G^i}{\partial \xi_j} + \Gamma_{lj}^i G^l\right) e_i\\
&= \left(F^3 \frac{\partial G^i}{\partial \xi_3} + P^i(F,G)\right) e_i,
\end{split}
\end{equation}
where $\Gamma_{ij}^k$'s are the Christoffel symbols of the second kind, and $P^i(F,G)$ collects all terms without $\xi_3$-derivative of $G$;\\
\item  Laplacian of a vector field:
\begin{equation}\label{laplacian}
\begin{split}
\Delta F &= E^{jk}\left\{\frac{\partial^2 F^i}{\partial \xi_j \partial \xi_k} - \Gamma_{jk}^h \frac{\partial F^i}{\partial \xi_h} + 2 \Gamma_{jh}^i \frac{\partial F^h}{\partial \xi_k} + F^h \frac{\partial}{\partial \xi_h} \Gamma_{jk}^i\right\} e_i\\
&= \left\{\frac{\partial^2 F^i}{\partial \xi_3^2} + S^i(F) + T^i(F)\right\} e_i,\\
\end{split}
\end{equation}
where $S^i(F)$ collects terms involving exactly one $\xi_3$-derivative, and $T^i(F)$ collects terms with no $\xi_3$-derivative (see \cite{Hirota1982});\\
\item Curl of a vector field:
\begin{equation}\label{curl}
(\operatorname{curl} F)^i = \frac{1}{\sqrt{V}} \left(\frac{\partial F_k}{\partial \xi_j} - \frac{\partial F_j}{\partial \xi_k}\right),
\end{equation}
where $(i,j,k)$ is a cyclic permutation of $(1,2,3)$.
\end{enumerate}

\bigskip

% section 3. the prandtl equation and boundary layer corrector
\section{The Prandtl Equation and Boundary Layer Corrector} \label{s:Prandtl}

The Prandtl equation formally governs the leading-order behavior in the asymptotic expansion of the NSE, and leads naturally to the construction of a boundary layer corrector.

Again formally, assuming $v^\epsilon \approx v^0 + \varphi^0$ and combining \eqref{NSE} with \eqref{EE}, we obtain the following system of equations for the corrector $\varphi^0$:
\begin{equation}
\begin{cases}
\partial_t \varphi^0 - \epsilon \Delta \varphi^0 + \mathcal{U} \cdot \nabla \varphi^0 + \varphi^0 \cdot \nabla \mathcal{U} + (v^0 \cdot \nabla \varphi^0 + \varphi^0 \cdot \nabla v^0 + \varphi^0 \cdot \nabla \varphi^0)\\
\quad + \nabla (p^\epsilon - p^0) \approx \epsilon \Delta (\mathcal{U} + v^0),\\
\operatorname{div} \varphi^0 = 0,\\
\gamma_- \varphi^0 = - \gamma_- v^0,\\
\gamma_+ \varphi^0 = 0,\\
\varphi^0|_{t=0} = 0.
\end{cases}
\end{equation}

Observe that $\mathcal{U}^3|_{\xi_3=0} = - \alpha^3$. Hence by the mean value formula, in each $H_{(k)}$
\begin{equation}\label{background flow BC}
\mathcal{U}^3 = - \alpha^3 + \frac{\partial \mathcal{U}^3}{\partial \xi_3} \Bigg|_{\xi_3=\zeta} \xi_3
\end{equation}
for some $\zeta=\zeta(k,\xi_3)$.

%\AnnaComment{At this point, we have not yet rigorously establish this below}
A detailed, but formal, asymptotic analysis, using \eqref{directional derivative}, \eqref{laplacian}, and \eqref{background flow BC}, shows that the boundary layer is expected to have thickness $O(\epsilon)$, and the dominant terms are given by:
\begin{align}
(- \epsilon \Delta \varphi_{(k)}^0 &+ \mathcal{U} \cdot \nabla \varphi_{(k)}^0)^i
= \left(- \epsilon \frac{\partial^2 \varphi_{(k)}^{0,i}}{\partial \xi_3^2} - \alpha^3 \frac{\partial \varphi_{(k)}^{0,i}}{\partial \xi_3}\right)  \nonumber \\
&+ \left(- \epsilon S^i(\varphi_{(k)}^0) + \frac{\partial \mathcal{U}^3}{\partial \xi_3} \Bigg|_{\xi=\zeta} \xi_3 \frac{\partial \varphi_{(k)}^{0,i}}{\partial \xi_3} + P^i(\mathcal{U}, \varphi_{(k)}^0)\right)
- \epsilon T^i(\varphi_{(k)}^0), 
\label{leading terms}
\end{align}
for $i=1,2,3$,  in each $H_{(k)}$. At leading order thus, the tangential components of the corrector formally satisfy the linear equation:
\begin{equation}\label{Prandtl}
- \epsilon \frac{\partial^2 \varphi_{(k)}^{0,i}}{\partial \xi_3^2} - \alpha^3 \frac{\partial \varphi_{(k)}^{0,i}}{\partial \xi_3} \approx 0,\qquad i=1,2,
\end{equation}
since the normal components of Euler and Navier-Stokes agree on the boundary. 
The error analysis carried out in Sections \ref{s:Sobolev} and \ref{s:correctors} will establish rigorously the validity of this expansion.

We first construct a local linear corrector, using \eqref{Prandtl}. To ensure incompressibility, we employ a vector potential. To this end, let $\{\eta_{(k)}\}_{k=1}^N$ be a partition of unity subordinated to $\{\Gamma_{(k)}\}_{k=1}^N$. We then set $\varphi_{(k)}^0 = \operatorname{curl} \psi_{(k)}^0$ in each $H_{(k)}$, where the potential function $\psi_{(k)}^0$ is defined by
\begin{equation}
\begin{cases}
\psi_{(k),i}^0 = (- 1)^{i-1} \epsilon \rho \eta_{(k)} \gamma_- (v^{0,3-i}) \frac{\sqrt{V_{(k)}}}{\alpha^3} \left(e^{-\frac{\alpha^3 \xi_3}{\epsilon}} - 1\right),\\
\psi_{(k),3}^0 = 0,
\end{cases}
\end{equation}
with $\rho = \rho(\xi_3)$ a smooth cutoff equal to $1$ on $[0,\frac{\delta}{2}]$ and $0$ on $[\frac{3\delta}{4},\delta]$. Applying \eqref{curl} gives the explicit formula:
\begin{equation}\label{BL formula}
\begin{cases}
\varphi_{(k)}^{0,i} = - \rho \eta_{(k)} \gamma_-(v^{0,i}) e^{-\frac{\alpha^3 \xi_3}{\epsilon}} + \epsilon \frac{\eta_{(k)} \gamma_-(v^{0,i})}{\alpha^3 \sqrt{V_{(k)}}} \frac{\partial}{\partial \xi_3} (\rho \sqrt{V_{(k)}}) \left(e^{-\frac{\alpha^3 \xi_3}{\epsilon}} - 1\right), i=1,2,\\
\begin{aligned}
\varphi_{(k)}^{0,3} = &- \epsilon \frac{\rho}{\sqrt{V_{(k)}}} \left[\frac{\partial}{\partial \xi_1} \left(\frac{\eta_{(k)} \gamma_-(v^{0,1}) \sqrt{V_{(k)}}}{\alpha^3}\right) + \frac{\partial}{\partial \xi_2} \left(\frac{\eta_{(k)} \gamma_-(v^{0,2}) \sqrt{V_{(k)}}}{\alpha^3}\right)\right] \left(e^{-\frac{\alpha^3 \xi_3}{\epsilon}} - 1\right)\\
&+ \xi_3 \frac{\eta_{(k)} \rho}{\alpha^3} \left[\gamma_-(v^{0,1})\frac{\partial \alpha^3}{\partial \xi_1} + \gamma_-(v^{0,2})\frac{\partial \alpha^3}{\partial \xi_2}\right] e^{-\frac{\alpha^3 \xi_3}{\epsilon}}.
\end{aligned}
\end{cases}
\end{equation}
%\AnnaComment{I am not particularly fond of this $C'$ notation as it can be confused with a constant. I would prefer to use something like $\sigma$ for example, though I understand you used $C'$ because it acts like a constant in $\epsilon$}
For notational convenience, we denote by $C$ a generic positive constant, by $\sigma$ a smooth function in $C^\infty_b(\Omega)$ with all derivatives uniformly bounded in $\epsilon$, and by $e.s.t.$ an exponentially small term along with its derivatives, i.e., bounded by $C e^{-\frac{C}{\epsilon^s}}$, for some $s>0$. Moreover, we introduce the notation $p_m(x,y)$ to represent a generic element of:
\begin{equation}
P_m = \left\{p(x,y): p(x,y) = \sum_{s\in\mathbb{Z}, t\in\mathbb{N}_0, m \leq s+t \leq M} \sigma x^s y^t, \mathbb{N}_0 \ni M \geq m\right\}, m \in \mathbb{N}_0.
\end{equation}
It is straightforward to verify that $P_{m+1} \subset P_m$ and
\begin{equation}
\frac{\partial}{\partial \xi_3} \left(p_{m+1}(\epsilon,\xi_3) e^{-\frac{r \alpha^3 \xi_3}{\epsilon}}\right) = p_m(\epsilon,\xi_3) e^{-\frac{r \alpha^3 \xi_3}{\epsilon}}, \forall r>0.
\end{equation}
This notation allows us to write \eqref{BL formula} in a more compact form:
\begin{equation}\label{BL formula 2}
\varphi_{(k)}^0 = \varphi_{(k),0}^0 + \varphi_{(k),1}^0 + \varphi_{(k),2}^0,\\
\end{equation}
where
\begin{equation}\label{BL formula 3}
\begin{cases}
\varphi_{(k),0}^{0,i} = - \rho \eta_{(k)} \gamma_-(v^{0,i}) e^{-\frac{\alpha^3 \xi_3}{\epsilon}}, \quad i=1,2,\\
\varphi_{(k),0}^{0,3} =0,\\
\varphi_{(k),1}^{0,i} = p_1(\epsilon,\xi_3) e^{-\frac{\alpha^3 \xi_3}{\epsilon}}, \quad i=1,2,3,\\
\varphi_{(k),2}^{0,i} = \sigma \epsilon, \quad  i=1,2,3,
\end{cases}
\end{equation}
with $p_1(x,y) = \sigma x + \sigma y$.

We note that $\varphi_{(k)}^0 \in C^\infty(H_{(k)})$ satisfies
\begin{equation}\label{excise BL}
- \epsilon \frac{\partial^2 \varphi_{(k)}^{0,i}}{\partial \xi_3^2} - \alpha^3 \frac{\partial \varphi_{(k)}^{0,i}}{\partial \xi_3} = p_0(\epsilon, \xi_3) e^{-\frac{\alpha^3 \xi_3}{\epsilon}} + \sigma \epsilon, \qquad i=1,2,3,
\end{equation}
where we have used the fact that any term of the form $\sigma \rho' e^{-\frac{C \xi_3}{\epsilon}}$ is an $e.s.t.$, and $p_0$ here takes the form $p_0(x,y) = \sigma$.

Furthermore, by construction and the compatibility condition \eqref{compatibility 2}, we have
\begin{equation}
\begin{cases}
\operatorname{div} \varphi_{(k)}^0 = 0,\\
\gamma_-(\varphi_{(k)}^0) = - \eta_{(k)} \gamma_-(v^0),\\
\varphi_{(k)}^0 \big|_{\partial H_{(k)} \backslash \Gamma_-} = 0,\\
\varphi_{(k)}^0 \big|_{t=0} = 0.
\end{cases}
\end{equation}
Therefore, we may extend $\varphi_{(k)}$ by zero smoothly to the full domain $\Omega$, and define the global boundary layer corrector by
\begin{equation} \label{eq:GlobalCorrector}
\varphi^0 = \sum_{i=1}^N \varphi_{(k)}^0.
\end{equation}
It follows that $\varphi \in C^\infty(\Omega)$ satisfies
\begin{equation}\label{BL}
\begin{cases}
\operatorname{div} \varphi^0 = 0,\\
\gamma_- \varphi^0 = - \gamma_- v^0,\\
\gamma_+ \varphi^0 = 0,\\
\varphi^0|_{t=0} = 0.
\end{cases}
\end{equation}

In the next section, we derive bounds in Sobolev spaces on the corrector, which in turn allows to estimate the difference between the NSE and EE in terms of $\epsilon$.

% section 4. sobolev estimates
\section{Sobolev Estimates} \label{s:Sobolev}

We first note that, for any $m \in \mathbb{N}_0$:
\begin{align}\label{BL control}
\left\|p_m(\epsilon,\xi_3) e^{-\frac{C \xi_3}{\epsilon}}\right\|_2 \leq C \epsilon^{\frac{1}{2}+m},\qquad 
\left\|p_m(\epsilon,\xi_3) e^{-\frac{C \xi_3}{\epsilon}}\right\|_\infty \leq C \epsilon^m. 
\end{align}

Setting $v^\epsilon = v^0 + \varphi^0 + w$ and substituting into \eqref{NSE}, \eqref{EE}, and \eqref{BL}, we derive a system for the remainder $w$:
\begin{equation}
\begin{cases}
\partial_t w - \epsilon \Delta w + (\mathcal{U} + v^\epsilon) \cdot \nabla w + \nabla (p^\epsilon - p^0) 
+ w \cdot \nabla (\mathcal{U} + v^0) + w \cdot \nabla \varphi^0= \epsilon \Delta (\mathcal{U} + v^0) \\
\qquad\qquad- \left[\partial_t \varphi^0 + (-\epsilon \Delta \varphi^0 + \mathcal{U} \cdot \nabla \varphi^0) + \varphi^0 \cdot \nabla (\mathcal{U} + v^0) + v^0 \cdot \nabla \varphi^0 + \varphi^0 \cdot \nabla \varphi^0\right] \\
\operatorname{div} w = 0,\\
\gamma w = 0, \label{remainder}\\
w|_{t=0} = 0.
\end{cases}
\end{equation}
Next, in order to estimate the right-hand side of the momentum equation for $w$, it is convenient to define
\begin{equation}\label{L_0}
\mathcal{L}^0 = \partial_t \varphi^0 + (-\epsilon \Delta \varphi^0 + \mathcal{U} \cdot \nabla \varphi^0) + \varphi^0 \cdot \nabla (\mathcal{U} + v^0) + v^0 \cdot \nabla \varphi^0 + \varphi^0 \cdot \nabla \varphi^0,
\end{equation}
%the collection of terms in the third line of \eqref{remainder}, which involve neither the remainder $w$ nor purely interior terms.\\
and rewrite it as  $\mathcal{L}^0 = \sum_{k=1}^N \mathcal{L}_{(k)}^0$, where for each $k$:
\begin{equation}
\mathcal{L}_{(k)}^0 := \partial_t \varphi_{(k)}^0 + (-\epsilon \Delta \varphi_{(k)}^0 + \mathcal{U} \cdot \nabla \varphi_{(k)}^0) + \varphi_{(k)}^0 \cdot \nabla (\mathcal{U} + v^0) + v^0 \cdot \nabla \varphi_{(k)}^0 + \varphi_{(k)}^0 \cdot \nabla \varphi^0.
\end{equation}

%\AnnaComment{I think that if we use the subscript $(k)$ on $\varphi$, we have to use it on $g_i$ as well}
We rewrite each term of $\mathcal{L}_{(k)}^0$ in  components as follows:
\begin{equation*}
\partial_t \varphi_{(k)}^0 = \sum_{i=1,2} \left(\sigma e^{-\frac{\alpha^3 \xi_3}{\epsilon}} + \sigma \epsilon\right) e_{(k),i} + \left(p_1(\epsilon,\xi_3) e^{-\frac{\alpha^3 \xi_3}{\epsilon}} + \sigma \epsilon\right) e_{(k),3}
\end{equation*}
(by \eqref{BL formula 2}, \eqref{BL formula 3});
\begin{align}
&- \epsilon \Delta \varphi_{(k)}^0 + \mathcal{U} \cdot \nabla \varphi_{(k)}^0\\
= &\left\{\left(- \epsilon \frac{\partial^2 \varphi_{(k)}^{0,i}}{\partial \xi_3^2} - \alpha^3 \frac{\partial \varphi_{(k)}^{0,i}}{\partial \xi_3}\right) + \left(- \epsilon S^i(\varphi_{(k)}^0) + \frac{\partial \mathcal{U}^3}{\partial \xi_3} \Bigg|_{\xi=\zeta} \xi_3 \frac{\partial \varphi_{(k)}^{0,i}}{\partial \xi_3} + P^i(\mathcal{U}, \varphi_{(k)}^0)\right) - \epsilon T^i(\varphi_{(k)}^0)\right\} e_{(k),i} \\
= &\sum_{i=1,2,3} \left(p_0(\epsilon,\xi_3) e^{-\frac{\alpha^3 \xi_3}{\epsilon}} + \sigma \epsilon\right) e_{(k),i} &&
\end{align}
(by \eqref{directional derivative}, \eqref{laplacian}, \eqref{leading terms}, \eqref{BL formula 2}, \eqref{BL formula 3}, \eqref{excise BL}, and $p_0$ here takes the form $p_0(x,y) = \sigma + \sigma \frac{y}{x}$);
\begin{equation*}
\begin{split}
\varphi_{(k)}^0 \cdot \nabla (\mathcal{U} + v^0) &= \varphi_{(k)}^{0,j} \left(\frac{\partial (\mathcal{U} + v^0)^i}{\partial \xi_j} + (\mathcal{U} + v^0)^l \Gamma_{lj}^i\right) e_{(k),i}\\
&= \sum_{i=1,2,3} \left(\sigma e^{-\frac{\alpha^3 \xi_3}{\epsilon}} + \sigma \epsilon\right) e_{(k),i}
\end{split}
\end{equation*}
(by \eqref{directional derivative}, \eqref{BL formula 2}, \eqref{BL formula 3});
\begin{equation}\label{L^0 1}
\begin{split}
v^0 \cdot \nabla \varphi_{(k)}^0 &= \left(\frac{\partial v^{0,3}}{\partial \xi_3} \bigg|_{\xi_3=\zeta} \xi_3 \frac{\partial \varphi_{(k)}^{0,i}}{\partial \xi_3} + P^i(v^0, \varphi_{(k)}^0)\right) e_{(k),i}\\
&= \sum_{i=1,2,3} \left(p_0(\epsilon,\xi_3) e^{-\frac{\alpha^3 \xi_3}{\epsilon}} + \sigma \epsilon\right) e_{(k),i}
\end{split}
\end{equation}
(by \eqref{directional derivative}, \eqref{BL formula 2}, \eqref{BL formula 3}, and $p_0$ here takes the form $p_0(x,y) = \sigma + \sigma \frac{y}{x}$).

The only term that requires some care is the quadratic term. Inserting  the decomposition \eqref{eq:GlobalCorrector} for $\varphi^0$, we distinguish two cases, when we multiply local correctors on the same chart or on overlapping charts (all other terms vanish).  In the first case, we have  $\varphi_{(k)}^0 \cdot \nabla \varphi_{(l)}^0$ with $k=l$. Then
\begin{equation}\label{L^0 2}
\begin{split}
\varphi_{(k)}^0 \cdot \nabla \varphi_{(k)}^0 &= \left(\varphi_{(k)}^{0,3} \frac{\partial \varphi_{(k)}^{0,i}}{\partial \xi_3} + P^i(\varphi_{(k)}^0, \varphi_{(k)}^0)\right) e_{(k),i}\\
&= \sum_{i=1,2,3} \left(\sigma e^{-\frac{\alpha^3 \xi_3}{\epsilon}} + p_0(\epsilon,\xi_3) e^{-\frac{2 \alpha^3 \xi_3}{\epsilon}} + \sigma \epsilon^2\right) e_{(k),i}
\end{split}
\end{equation}
(by \eqref{directional derivative}, \eqref{BL formula 2}, \eqref{BL formula 3}, and $p_0$ here takes the form $p_0(x,y) = \sigma + \sigma \frac{y}{x}$).
In the second case, we have  $\varphi_{(k)}^0 \cdot \nabla \varphi_{(l)}^0$ with $k\ne l$ and $H_{(k)} \cap H_{(l)} \ne \emptyset$. Since (contravariant) vectors are transformed from one chart to the other via a matrix of the form:
\begin{equation*}
T = 
\begin{pmatrix}
t_{11} & t_{12} & 0 \\
t_{21} & t_{22} & 0 \\
0 & 0 & 1
\end{pmatrix},
\end{equation*}
the structure of the transformed corrector remains the same and we can perform a similar analysis on $\varphi_{(k)}^0 \cdot \nabla \varphi_{(l)}^0$. As a consequence, we can write
\begin{equation}\label{linear terms}
\mathcal{L}_{(k)}^0 = \sum_{i=1,2,3} \left(p_0(\epsilon,\xi_3) e^{-\frac{\alpha^3 \xi_3}{\epsilon}} + p_0(\epsilon,\xi_3) e^{-\frac{2 \alpha^3 \xi_3}{\epsilon}} + \sigma \epsilon\right) e_{(k),i},
\end{equation}
with $p_0(x,y) = \sigma + \sigma \frac{y}{x}$.

We can now perform energy estimates on \eqref{remainder}. We first observe that
\begin{equation*}
\Big|\int - \epsilon \Delta (\mathcal{U} + v^0) \cdot w\Big| \leq C \epsilon^2 + C \|w\|_2^2,
\end{equation*}
and
\begin{equation*}
\begin{split}
\Big|\int \mathcal{L}_{(k)}^0 \cdot w\Big| &= \Big|\int \mathcal{L}_{(k)}^{0,i} w_i\Big|
= \sum_{i=1,2,3} \Big|\int \left(p_0(\epsilon,\xi_3) e^{-\frac{\alpha^3 \xi_3}{\epsilon}} + p_0(\epsilon,\xi_3) e^{-\frac{2 \alpha^3 \xi_3}{\epsilon}} + \sigma \epsilon\right) w_i\Big|\\
&\leq C \left(\left\|p_1(\epsilon,\xi_3) e^{-\frac{\alpha_0 \xi_3}{\epsilon}}\right\|_2 + \left\|p_1(\epsilon,\xi_3) e^{-\frac{2 \alpha_0 \xi_3}{\epsilon}}\right\|_2\right) \left\|\frac{w}{\xi_3}\right\|_2 + C \epsilon \|w\|_2\\
&\leq C \epsilon^{3/2} \|\nabla w\|_2^2 + C \epsilon \|w\|_2\\
&\leq C \epsilon^2 + \frac{\epsilon}{4N} \|\nabla w\|_2^2 + C \|w\|_2^2,
\end{split}
\end{equation*}
where we used \eqref{BL control} and Hardy's inequality.
On the other hand, 
\begin{equation}\label{nonlinear 1}
\Big|\int (w \cdot \nabla (\mathcal{U} + v^0)) \cdot w\Big| \leq C \|w\|_2^2,
\end{equation}
and
\begin{equation}\label{nonlinear 2}
\begin{split}
\Big|\int (w \cdot \nabla \varphi_{(k)}^0) \cdot w\Big| \leq &\Big|\int w^3 \frac{\partial \varphi_{(k),0}^{0,i}}{\partial \xi_3} w_i\Big| + \Big|\int w^3 \frac{\partial \left(\varphi_{(k),1}^{0,i} + \varphi_{(k),2}^{0,i}\right)}{\partial \xi_3} w_i + P^i(w, \varphi_{(k)}^0) w_i\Big|\\
\leq &\Big|\int \frac{\partial w^3}{\partial \xi_3} \varphi_{(k),0}^{0,i} w_i\Big| + \Big|\int w^3 \varphi_{(k),0}^{0,i} \frac{\partial w_i}{\partial \xi_3}\Big|\\
&+ \Big|\int w^3 \frac{\partial \left(\varphi_{(k),1}^{0,i} + \varphi_{(k),2}^{0,i}\right)}{\partial \xi_3} w_i + P^i(w, \varphi_{(k)}^0) w_i\Big|\\
\leq &\Big|\int \frac{\partial w^3}{\partial \xi_3} \varphi_{(k),0}^{0,i} w_i\Big| + \Big|\int w^3 \varphi_{(k),0}^{0,i} \frac{\partial w^i}{\partial \xi_3}\Big|\\
&+ \sum _{i,j=1,2,3} \Big|\int w^j \left(p_0(\epsilon,\xi_3) e^{- \frac{\alpha_0 \xi_3}{\epsilon}} + \sigma \epsilon\right) w_i\Big|\\
\leq &\left\|\frac{\partial w^3}{\partial \xi_3}\right\|_2 \left\|\xi_3 \varphi_{(k),0}^{0,i}\right\|_\infty \left\|\frac{w_i}{\xi_3}\right\|_2 + \left\|\frac{w^3}{\xi_3}\right\|_2 \left\|\xi_3 \varphi_{(k),0}^{0,i}\right\|_\infty \left\|\frac{\partial w_i}{\partial \xi_3}\right\|_2\\
&+ \left\|p_0(\epsilon,\xi_3) e^{- \frac{\alpha_0 \xi_3}{\epsilon}} + \sigma \epsilon\right\|_\infty \|w\|_2^2\\
\leq &\frac{\epsilon}{4N} \|\nabla w\|_2^2 + C \|w\|_2^2,
\end{split}
\end{equation}
where in the last step we used the assumption
\begin{equation}
\|\gamma_-(v^0)\|_\infty \leq \frac{e \alpha_0}{8N},
\end{equation}
which is valid for sufficiently small $T>0$, ensuring that
\begin{equation}\label{small time}
\begin{split}
\|\xi_3 \varphi_{(k),0}^{0,i}\|_\infty &= \left\|\xi_3 \rho \eta_{(k)} \gamma_-(v^{0,i}) e^{-\frac{\alpha^3 \xi_3}{\epsilon}}\right\|_\infty\\
&\leq \left\|\xi_3 e^{-\frac{\alpha_0 \xi_3}{\epsilon}}\right\|_\infty \|\gamma_-(v^{0,i})\|_\infty\\
&\leq \frac{\epsilon}{e \alpha_0} \|\gamma_-(v^0)\|_\infty
\leq \frac{\epsilon}{8N}.
\end{split}
\end{equation}
From the estimates above,  multiplying the first equation in $\eqref{remainder}$ by $w$ and integrating over $\Omega$ yields
\begin{equation*}
\frac{d}{dt} \|w\|_2^2 + \epsilon \|\nabla w\|_2^2 \leq C \epsilon^2 + C \|w\|_2^2,
\end{equation*}
from which it follows by Gr\" onwall's Lemma  that
\begin{equation}\label{cvg rate}
%\begin{cases}
\|w\|_{L^\infty L^2} \leq \epsilon,\quad
\|w\|_{L^2 H^1} \leq \epsilon^{1/2}.
%\end{cases}
\end{equation}
We note that, as typical of situations where the boundary layer is linear, there is no need to correct the pressure at leading order (cf. \cite{Gie2010,GieKelliherNavier2012,GieEtAl2018}).

% section 5. the higher order boundary layer correctors
\section{Higher-Order Boundary Layer Correctors} \label{s:correctors}

In this section, we construct a refined boundary layer corrector by performing a higher-order asymptotic expansion of the velocity and pressure. 
Using higher-order correctors leads to an improvement of the convergence rates established in \eqref{cvg rate}, raising them by one order. Specifically, we write
\begin{equation}\label{higher order expansion}
\begin{cases}
v^\epsilon \approx (v^0 + \varphi^0) + \epsilon (v^1 + \varphi^1),\\
p^\epsilon \approx p^0 + \epsilon (p^1 + q^1),
\end{cases}
\end{equation}
where $v^1$, $\varphi^1$, $p^1$, and $q^1$ are the first-order correctors, to be determined later. We introduce two correctors to separate the behavior in the bulk, captured by $v^1$ and $p^1$, respectively, from the behavior in the boundary layer, captured by $\varphi^1$ and $q^1$, respectively. These last two correctors will be constructed locally using the curvilinear coordinate systems and the partition of unity already introduced.

From the analysis in Section \ref{s:Sobolev}, we recall that $\mathcal{L}_{(k)}^0$ can be written as
\begin{equation}
\mathcal{L}_{(k)}^0 = - \left(\sigma_{(k),0}^i \epsilon + \sigma_{(k),1}^i e^{-\frac{\alpha^3 \xi_3}{\epsilon}} + \sigma_{(k),2}^i \frac{\xi_3}{\epsilon} e^{-\frac{\alpha^3 \xi_3}{\epsilon}} + \sigma_{(k),3}^i e^{-\frac{2 \alpha^3 \xi_3}{\epsilon}} + \sigma_{(k),4}^i \frac{\xi_3}{\epsilon} e^{-\frac{2 \alpha^3 \xi_3}{\epsilon}}\right) e_{(k),i},
\end{equation}
where each coefficient $\sigma_{(k),j}^i$ belongs to $C^\infty_b(\Omega)$ with all derivatives uniformly bounded in $\epsilon$, while $e_{(k),i}$ are the elements of the frame induced by the curvilinear coordinate system. By construction and the compatibility conditions \eqref{compatibility 2} and \eqref{compatibility 4}, $\mathcal{L}_{(k)}^0$ vanishes outside $H_{(k)}$ and at $t = 0$. Consequently, each $\sigma_{(k),j}^i$ shares these properties.

For each $i$ and $j \ne 0$, we apply again the mean value theorem to conclude that
\begin{equation}
\sigma_{(k),j}^i = \sigma_{(k),j}^{'i} + \frac{\partial \sigma_{(k),j}^{i}}{\partial \xi_3} \Bigg|_{\xi_3=\zeta} \xi_3,
\end{equation}
where $\sigma_{(k),j}^{'i} := \gamma_- \sigma_{(k),j}^i$, for some $\zeta=\zeta(k,\xi_3)$. This decomposition ensures that $\sigma_{(k),j}^{'i}$ is independent of $\xi_3$ and vanishes outside $\Gamma_{(k)}$ and at $t = 0$.
We can thus rewrite
\begin{equation}
\mathcal{L}_{(k)}^0 = - \epsilon \mathcal{Q}_{(k)} - \mathcal{M}_{(k)} + \mathcal{R}_{(k)},
\end{equation}
with
\begin{equation}\label{R_k}
\begin{split}
\mathcal{Q}_{(k)} = &\sigma_{(k),0}^i e_{(k),i},\\
\mathcal{M}_{(k)} = &\left(\sigma_{(k),1}^{'i} e^{-\frac{\alpha^3 \xi_3}{\epsilon}} + \sigma_{(k),2}^{'i} \frac{\xi_3}{\epsilon} e^{-\frac{\alpha^3 \xi_3}{\epsilon}} + \sigma_{(k),3}^{'i} e^{-\frac{2 \alpha^3 \xi_3}{\epsilon}} + \sigma_{(k),4}^{'i} \frac{\xi_3}{\epsilon} e^{-\frac{2 \alpha^3 \xi_3}{\epsilon}}\right) e_{(k),i},\\
\mathcal{R}_{(k)} = &- \left(\frac{\partial \sigma_{(k),1}^i}{\partial \xi_3} \Bigg|_{\xi_3=\zeta} \xi_3 e^{-\frac{\alpha^3 \xi_3}{\epsilon}} + \frac{\partial \sigma_{(k),2}^i}{\partial \xi_3} \Bigg|_{\xi_3=\zeta} \frac{\xi_3^2}{\epsilon} e^{-\frac{\alpha^3 \xi_3}{\epsilon}}\right.\\
&\left.+ \frac{\partial \sigma_{(k),3}^i}{\partial \xi_3} \Bigg|_{\xi_3=\zeta} \xi_3 e^{-\frac{2 \alpha^3 \xi_3}{\epsilon}} + \frac{\partial \sigma_{(k),4}^i}{\partial \xi_3} \Bigg|_{\xi_3=\zeta} \frac{\xi_3^2}{\epsilon} e^{-\frac{2 \alpha^3 \xi_3}{\epsilon}} \right) e_{(k),i}\\
= &\sum_{i=1,2,3} \left(p_1(\epsilon,\xi_3) e^{-\frac{\alpha^3 \xi_3}{\epsilon}} + p_1(\epsilon,\xi_3) e^{-\frac{2 \alpha^3 \xi_3}{\epsilon}}\right) e_{(k),i}.
\end{split}
\end{equation}

Writing as a short-hand notation
\begin{equation}
\begin{split}
\mathcal{Q} = \sum_{k=1}^N \mathcal{Q}_{(k)}, \quad \mathcal{M} = \sum_{k=1}^N \mathcal{M}_{(k)}, \quad \mathcal{R} = \sum_{k=1}^N \mathcal{R}_{(k)},
\end{split}
\end{equation}
we have that
\begin{equation}
\begin{split}
\mathcal{L}^0 = \sum_{k=1}^N \mathcal{L}_{(k)} = - \epsilon \mathcal{Q} - \mathcal{M} + \mathcal{R}.
\end{split}
\end{equation}

% 5.1. the higher order bulk equation
\subsection{The higher-order corrector equations}

Introducing the expansion in \eqref{higher order expansion} into the Navier-Stokes equation \eqref{NSE}, we obtain the following system of equations:
\begin{align}
&\partial_t ((v^0 + \varphi^0) + \epsilon (v^1 + \varphi^1)) - \epsilon \Delta ((v^0 + \varphi^0) + \epsilon (v^1 + \varphi^1)) \nonumber \\
&+ \mathcal{U} \cdot \nabla ((v^0 + \varphi^0) + \epsilon (v^1 + \varphi^1)) + ((v^0 + \varphi^0) + \epsilon (v^1 + \varphi^1)) \cdot \nabla \mathcal{U} \nonumber \\
&+ ((v^0 + \varphi^0) + \epsilon (v^1 + \varphi^1)) \cdot \nabla ((v^0 + \varphi^0) + \epsilon (v^1 + \varphi^1)) + \nabla (p^0 + \epsilon (p^1 + q^1))\nonumber \\
&\approx f + \epsilon \Delta \mathcal{U}. \label{1st approx NSE}
\end{align}

Subtracting the Euler equation \eqref{EE}, satisfied by the leading order term $v^0$, from \eqref{1st approx NSE} and incorporating the expression \eqref{L_0} for $\mathcal{L}_0$, we formally obtain the system  for $v^1$:
\begin{equation}\label{bulk 1}
\begin{cases}
\partial_t v^1 + (\mathcal{U} + v^0) \cdot \nabla v^1 + v^1 \cdot \nabla (\mathcal{U} + v^0) + \nabla p^1 = \mathcal{Q} + \Delta (\mathcal{U} + v_0),\\
\operatorname{div} v^1 = 0,\\
\gamma_+ v^1 = 0,\\
\gamma_- v^1 \cdot n = 0,\\
v^1|_{t=0} = 0,
\end{cases}
\end{equation}
which admits a unique smooth solution. Above we have separated out all terms containing $\varphi^1$ and $q^1$, which will be treated next.

% 5.2. the higher order Prandtl equation and boundary layer correctors
\subsection{The higher order Prandtl equation and boundary layer correctors}

Using \eqref{bulk 1} into \eqref{1st approx NSE} and utilizing \eqref{directional derivative}, \eqref{laplacian}, \eqref{background flow BC}, \eqref{L_0}, and \eqref{R_k}, allow to identify the dominant terms in the remainder.
Localizing to each coordinate patch $H_{(k)}$, and denoting with $\varphi_{(k)}^1$ the local corrector, we have:
\begin{align}
&\epsilon \left(- \epsilon \Delta \varphi_{(k)}^1 + \mathcal{U} \cdot \nabla \varphi_{(k)}^1 - \frac{1}{\epsilon} \mathcal{M}_{(k)}\right)^i
= \epsilon \left\{\left(- \epsilon \frac{\partial^2 \varphi_{(k)}^{1,i}}{\partial \xi_3^2} - \alpha^3 \frac{\partial \varphi_{(k)}^{1,i}}{\partial \xi_3} - \frac{1}{\epsilon} \mathcal{M}_{(k)}^i\right)\right. \label{excise higher order BL} \\
&\left.+ \left(- \epsilon S^i(\varphi_{(k)}^1) + \frac{\partial \mathcal{U}^3}{\partial \xi_3} \Big|_{\xi_3 = \zeta} \xi_3 \frac{\partial \varphi_{(k)}^{1,i}}{\partial \xi_3} + P^i(\mathcal{U}, \varphi_{(k)}^1)\right) - \epsilon T^i(\varphi_{(k)}^1)\right\}, \qquad i=1,2,3. \nonumber
\end{align}

Collecting the leading-order terms gives the following linear ODE for $\varphi_{(k)}^1$:
\begin{align}
&- \epsilon \frac{\partial^2 \varphi_{(k)}^{1,i}}{\partial \xi_3^2} - \alpha^3 \frac{\partial \varphi_{(k)}^{1,i}}{\partial \xi_3} \approx \frac{1}{\epsilon} \mathcal{M}_{(k)}^i \label{higher order Prandtl} \\
&= \frac{1}{\epsilon} \left(\sigma_{(k),1}^{'i} e^{-\frac{\alpha^3 \xi_3}{\epsilon}} + \sigma_{(k),2}^{'i} \frac{\xi_3}{\epsilon} e^{-\frac{\alpha^3 \xi_3}{\epsilon}} + \sigma_{(k),3}^{'i} e^{-\frac{2 \alpha^3 \xi_3}{\epsilon}} + \sigma_{(k),4}^{'i} \frac{\xi_3}{\epsilon} e^{-\frac{2 \alpha^3 \xi_3}{\epsilon}}\right), \qquad i=1,2. \nonumber
\end{align}
Solving this ODE gives
\begin{equation}
\begin{split}
\varphi_{(k)}^{1,i} \approx &\left(\frac{\sigma_{(k),3}^{'i}}{2(\alpha^3)^2} + \frac{3 \sigma_{(k),4}^{'i}}{4 (\alpha^3)^3} - \eta_{(k)} \gamma_-(v^{1,i})\right) e^{-\frac{\alpha^3 \xi_3}{\epsilon}}\\
&+ \left[\left(\frac{\sigma_{(k),1}^{'i}}{\alpha^3} + \frac{\sigma_{(k),2}^{'i}}{(\alpha^3)^2}\right) \left(\frac{\xi_3}{\epsilon}\right) + \frac{\sigma_{(k),2}^{'i}}{2 \alpha^3} \left(\frac{\xi_3}{\epsilon}\right)^2\right] e^{-\frac{\alpha^3 \xi_3}{\epsilon}}\\
&- \left[\left(\frac{\sigma_{(k),3}^{'i}}{2 (\alpha^3)^2} + \frac{3 \sigma_{(k),4}^{'i}}{4 (\alpha^3)^3}\right) + \frac{\sigma_{(k),4}^{'i}}{2 (\alpha^3)^2} \left(\frac{\xi_3}{\epsilon}\right)\right] e^{-\frac{2 \alpha^3 \xi_3}{\epsilon}}, \qquad i=1,2.
\end{split}
\end{equation}

We can compute the third component by enforcing the divergence-free condition. To this end, it is again conveneint to introduce a vector potential $\psi_{(k)}^1,$ which satisfies $\varphi_{(k)}^1 = \operatorname{curl} \psi_{(k)}^1$
and is given by:
\begin{align}
\psi_{(k),i}^{1} = &(-1)^{i-1} \rho \sqrt{V_{(k)}} \left\{-\frac{\epsilon}{\alpha^3} \left(\frac{\sigma_{(k),3}^{'i}}{2(\alpha^3)^2} + \frac{3 \sigma_{(k),4}^{'i}}{4 (\alpha^3)^3} - \eta_{(k)} \gamma_-(v^{1,3-i})\right) e^{-\frac{\alpha^3 \xi_3}{\epsilon}}\right.\\
&+ \left[\left(- \frac{\epsilon}{(\alpha^3)^2} -\frac{\xi_3}{\alpha^3}\right) \left(\frac{\sigma_{(k),1}^{'i}}{\alpha^3} + \frac{\sigma_{(k),2}^{'i}}{(\alpha^3)^2}\right) + \left(- \frac{2 \epsilon}{(\alpha^3)^3} -\frac{\xi_3^2}{\epsilon \alpha^3} - \frac{2 \xi_3}{(\alpha^3)^2}\right) \frac{\sigma_{(k),2}^{'i}}{2 \alpha^3}\right] e^{-\frac{\alpha^3 \xi_3}{\epsilon}}\\
&+ \left[\frac{\epsilon}{2 \alpha^3} \left(\frac{\sigma_{(k),3}^{'i}}{2 (\alpha^3)^2} + \frac{3 \sigma_{(k),4}^{'i}}{4 (\alpha^3)^3}\right) + \left(\frac{\epsilon}{4 (\alpha^3)^2} + \frac{\xi_3}{2 \alpha^3}\right) \frac{\sigma_{(k),4}^{'i}}{2 (\alpha^3)^2}\right] e^{-\frac{2 \alpha^3 \xi_3}{\epsilon}}\\
&\left.\epsilon \left(- \frac{\eta_{(k)} \gamma_-(v^{1,3-i})}{\alpha^3} + \frac{\sigma_{(k),1}^{'i}}{(\alpha^3)^3} + \frac{2 \sigma_{(k),2}^{'i}}{(\alpha^3)^4} + \frac{\sigma_{(k),3}^{'i}}{4 (\alpha^3)^3} + \frac{\sigma_{(k),4}^{'i}}{4 (\alpha^3)^4}\right)\right\}, \nonumber \\
\psi_{(k),3}^{1} = &0. && \label{higher order potential}
\end{align}

We can rewrite the expression above in compact form as  
\begin{equation*}
\psi_{(k),i}^1 = p_1(\epsilon,\xi_3) e^{-\frac{\alpha^3 \xi_3}{\epsilon}} + p_1(\epsilon,\xi_3) e^{-\frac{2 \alpha^3 \xi_3}{\epsilon}} + \sigma \epsilon, i=1,2,\\
\psi_{(k),3}^1 = 0.
\end{equation*}
Applying the curl yields
\begin{align}
\varphi_{(k)}^{1,i} = &\rho \left\{\left(\frac{\sigma_{(k),3}^{'i}}{2(\alpha^3)^2} + \frac{3 \sigma_{(k),4}^{'i}}{4 (\alpha^3)^3} - \eta_{(k)} \gamma_-(v^{1,i})\right) e^{-\frac{\alpha^3 \xi_3}{\epsilon}}\right. \nonumber \\
&+ \left[\left(\frac{\sigma_{(k),1}^{'i}}{\alpha^3} + \frac{\sigma_{(k),2}^{'i}}{(\alpha^3)^2}\right) \left(\frac{\xi_3}{\epsilon}\right) + \frac{\sigma_{(k),2}^{'i}}{2 \alpha^3} \left(\frac{\xi_3}{\epsilon}\right)^2\right] e^{-\frac{\alpha^3 \xi_3}{\epsilon}}  \nonumber \\
&\left.- \left[\left(\frac{\sigma_{(k),3}^{'i}}{2 (\alpha^3)^2} + \frac{3 \sigma_{(k),4}^{'i}}{4 (\alpha^3)^3}\right) + \frac{\sigma_{(k),4}^{'i}}{2 (\alpha^3)^2} \left(\frac{\xi_3}{\epsilon}\right)\right] e^{-\frac{2 \alpha^3 \xi_3}{\epsilon}}\right\} \nonumber \\
&+ (-1)^{i-1} \frac{\partial}{\partial \xi_3} (\rho \sqrt{V_{(k)}}) \frac{1}{\rho \sqrt{V_{(k)}}} \psi_{(k),3-i}^1, \qquad i = 1,2, \nonumber \\
\varphi_{(k)}^{1,3} = &\frac{1}{\sqrt{V_{(k)}}} \left(\frac{\partial \psi_{(k),2}^1}{\partial \xi_1} - \frac{\partial \psi_{(k),1}^1}{\partial \xi_2}\right),
\label{higher order BL formula}
\end{align}
which can be rewritten in more compact form as
\begin{align}
\varphi_{(k)}^{1,i} = &\rho \left\{\left(\frac{\sigma_{(k),3}^{'i}}{2(\alpha^3)^2} + \frac{3 \sigma_{(k),4}^{'i}}{4 (\alpha^3)^3} - \eta_{(k)} \gamma_-(v^{1,i})\right) e^{-\frac{\alpha^3 \xi_3}{\epsilon}}\right. \nonumber\\
&+ \left[\left(\frac{\sigma_{(k),1}^{'i}}{\alpha^3} + \frac{\sigma_{(k),2}^{'i}}{(\alpha^3)^2}\right) \left(\frac{\xi_3}{\epsilon}\right) + \frac{\sigma_{(k),2}^{'i}}{2 \alpha^3} \left(\frac{\xi_3}{\epsilon}\right)^2\right] e^{-\frac{\alpha^3 \xi_3}{\epsilon}} \nonumber \\
&\left.- \left[\left(\frac{\sigma_{(k),3}^{'i}}{2 (\alpha^3)^2} + \frac{3 \sigma_{(k),4}^{'i}}{4 (\alpha^3)^3}\right) + \frac{\sigma_{(k),4}^{'i}}{2 (\alpha^3)^2} \left(\frac{\xi_3}{\epsilon}\right)\right] e^{-\frac{2 \alpha^3 \xi_3}{\epsilon}}\right\} \nonumber\\
&+ p_1(\epsilon,\xi_3) e^{-\frac{\alpha^3 \xi_3}{\epsilon}} + p_1(\epsilon,\xi_3) e^{-\frac{2 \alpha^3 \xi_3}{\epsilon}} + \sigma \epsilon, \qquad i=1,2, \nonumber\\
\varphi_{(k)}^{1,3} = &p_1(\epsilon,\xi_3) e^{-\frac{\alpha^3 \xi_3}{\epsilon}} + p_1(\epsilon,\xi_3) e^{-\frac{2 \alpha^3 \xi_3}{\epsilon}} + \sigma \epsilon. \label{higher order BL}
\end{align}

We observe that $\varphi_{(k)}^1 \in C^\infty(H_{(k)})$ satisfies
\begin{equation}\label{excise higher order BL 2}
\begin{cases}
- \epsilon \frac{\partial^2 \varphi_{(k)}^{1,i}}{\partial \xi_3^2} - \alpha^3 \frac{\partial \varphi_{(k)}^{1,i}}{\partial \xi_3} = \frac{1}{\epsilon} \mathcal{M}_{(k)}^i + p_0(\epsilon,\xi_3) e^{-\frac{\alpha^3 \xi_3}{\epsilon}} + p_0(\epsilon,\xi_3) e^{-\frac{2 \alpha^3 \xi_3}{\epsilon}} + \sigma \epsilon, i=1,2,\\
- \epsilon \frac{\partial^2 \varphi_{(k)}^{1,3}}{\partial \xi_3^2} - \alpha^3 \frac{\partial \varphi_{(k)}^{1,3}}{\partial \xi_3} = p_0(\epsilon,\xi_3) e^{-\frac{\alpha^3 \xi_3}{\epsilon}} + p_0(\epsilon,\xi_3) e^{-\frac{2 \alpha^3 \xi_3}{\epsilon}} + \sigma \epsilon,
\end{cases}
\end{equation}
and
\begin{equation*}
\begin{cases}
\operatorname{div} \varphi_{(k)}^1 = 0,\\
\gamma_-\varphi_{(k)}^1 = - \eta_{(k)} \gamma_-v^1,\\
\varphi_{(k)}^1\big|_{\partial H_{(k)} \backslash \Gamma_-} = 0,\\
\varphi_{(k)}^1\big|_{t=0} = 0,
\end{cases}
\end{equation*}
so we may extend $\varphi_{(k)}^1$ by zero smoothly to the full domain $\Omega$ and define the global corrector
\begin{equation*}
\varphi^1 = \sum_{i=1}^N \varphi_{(k)}^1.
\end{equation*}
This corrector has the desired properties that 
\begin{equation}\label{BL2}
\begin{cases}
\operatorname{div} \varphi^1 = 0,\\
\gamma_- \varphi^1 = - \gamma_- v^1,\\
\gamma_+ \varphi^1 = 0,\\
\varphi^1|_{t=0} = 0.
\end{cases}
\end{equation}

Finally, since \eqref{higher order Prandtl} holds only for \(i=1,2\), we introduce a local pressure corrector \(q_{(k)}^1\) to balance the $i=3$ component, satisfying
\begin{align}
(\nabla q_{(k)}^1)^3 &= \frac{\partial q_{(k)}^1}{\partial \xi_3} = \frac{1}{\epsilon} \mathcal{M}_{(k)}^3 
\nonumber \\
&= \frac{1}{\epsilon} \left(\sigma_{(k),1}^{'3} e^{-\frac{\alpha^3 \xi_3}{\epsilon}} + \sigma_{(k),2}^{'3} \frac{\xi_3}{\epsilon} e^{-\frac{\alpha^3 \xi_3}{\epsilon}} + \sigma_{(k),3}^{'3} e^{-\frac{2 \alpha^3 \xi_3}{\epsilon}} + \sigma_{(k),4}^{'3} \frac{\xi_3}{\epsilon} e^{-\frac{2 \alpha^3 \xi_3}{\epsilon}}\right). \nonumber
\end{align}

Integrating yields
\begin{align}
q_{(k)}^1 = &- \frac{1}{\alpha^3} C_{(k),1}^3 e^{-\frac{\alpha^3 \xi_3}{\epsilon}} - C_{(k),2}^3 \left(\frac{1}{(\alpha^3)^2} + \frac{1}{\alpha^3} \frac{\xi_3}{\epsilon}\right) e^{-\frac{\alpha^3 \xi_3}{\epsilon}} \nonumber \\ 
&- \frac{1}{2 \alpha^3} C_{(k),3}^3 e^{-\frac{2 \alpha^3 \xi_3}{\epsilon}} - C_{(k),4}^3 \left(\frac{1}{4 (\alpha^3)^2} + \frac{1}{2 \alpha^3} \frac{\xi_3}{\epsilon}\right) e^{-\frac{2 \alpha^3 \xi_3}{\epsilon}}, \nonumber
\end{align}
or in a compact form,
\begin{equation*}
q_{(k)}^1 = p_0(\epsilon,\xi_3) e^{-\frac{\alpha^3 \xi_3}{\epsilon}} + p_0(\epsilon,\xi_3) e^{-\frac{2 \alpha^3 \xi_3}{\epsilon}}.
\end{equation*}

By \eqref{gradient} we then obtain:
\begin{equation}\label{grad pressure}
\nabla q_{(k)}^1 = \sum_{i=1,2} \left(p_0(\epsilon,\xi_3) e^{-\frac{\alpha^3 \xi_3}{\epsilon}} + p_0(\epsilon,\xi_3) e^{-\frac{2 \alpha^3 \xi_3}{\epsilon}}\right) e_{(k),i} + \frac{1}{\epsilon} \mathcal{M}_{(k)}^3 e_{(k),3}.
\end{equation}

% 5.3. higher order sobolev estimates
\subsection{Higher order Sobolev estimates}

We conclude this section by deriving error estimates for the first-order expansion that we have formally obtained. These estimates rigorously justify the asymptotic expension to order one.

We hence set $v^\epsilon = (v^0 + \varphi^0) + \epsilon (v^1 + \varphi^1) + w$, $p^\epsilon = p^0 + \epsilon (p^1 + q^1) + r$, and substitute this expression  into \eqref{NSE}, where $w$ and $r$ are the remainders for the velocity and pressure, repsectively. Subtracting \eqref{EE}, \eqref{BL}, \eqref{bulk 1}, and \eqref{BL2}, we obtain the system for $w$:
\begin{equation}\label{higher order remainder}
\begin{cases}
\partial_t w - \epsilon \Delta w + (\mathcal{U} + v^\epsilon) \cdot \nabla w + \nabla r\\
\quad + \epsilon^2 (- \Delta v^1 + v^1 \cdot \nabla v^1)\\
\quad + \epsilon \left\{\left(-\epsilon \Delta \varphi^1 + \mathcal{U} \cdot \nabla \varphi^1 - \frac{1}{\epsilon} \mathcal{M} + \nabla q^1\right)\right.\\
\quad + \left(\frac{1}{\epsilon} \mathcal{R} + \partial_t \varphi^1 + \varphi^1 \cdot \nabla (\mathcal{U} + v^0 + \epsilon v^1) + \varphi^0 \cdot \nabla v^1\right)\\
\left.\quad + (v^1 \cdot \nabla \varphi^0 + v^0 \cdot \nabla \varphi^1) + (\varphi^1 \cdot \nabla \varphi^0 + \varphi^0 \cdot \nabla \varphi^1) + \epsilon (v^1 \cdot \nabla \varphi^1 + \varphi^1 \cdot \nabla \varphi^1)\right\}\\
\quad + w \cdot \nabla (\mathcal{U} + v^0 + \varphi^0) + \epsilon w \cdot \nabla (v^1 + \varphi^1) =0,\\
\operatorname{div} w = 0,\\
\gamma w = 0,\\
w|_{t=0} = 0.
\end{cases}
\end{equation}

We now define:
\begin{align}
\mathcal{L}^1 = &\epsilon \left\{\left(-\epsilon \Delta \varphi^1 + \mathcal{U} \cdot \nabla \varphi^1 - \frac{1}{\epsilon} \mathcal{M} + \nabla q^1\right)\right. \nonumber \\
&+ \left(\frac{1}{\epsilon} \mathcal{R} + \partial_t \varphi^1 + \varphi^1 \cdot \nabla (\mathcal{U} + v^0 + \epsilon v^1) + \varphi^0 \cdot \nabla v^1\right) \nonumber \\
&\left.+ (v^1 \cdot \nabla \varphi^0 + v^0 \cdot \nabla \varphi^1) + (\varphi^1 \cdot \nabla \varphi^0 + \varphi^0 \cdot \nabla \varphi^1) + \epsilon (v^1 \cdot \nabla \varphi^1 + \varphi^1 \cdot \nabla \varphi^1)\right\}. \nonumber
\end{align}
$\mathcal{L}^1$ is the collection of all terms from line 3-5 of \eqref{higher order remainder}, which involves neither $w$ nor purely interior contributions. For each $k$, we introduce
\begin{align}
\mathcal{L}_{(k)}^1 := &\epsilon \left\{\left(-\epsilon \Delta \varphi_{(k)}^1 + \mathcal{U} \cdot \nabla \varphi_{(k)}^1 - \frac{1}{\epsilon} \mathcal{M}_{(k)} + \nabla q_{(k)}^1\right)\right. \nonumber \\
&+ \left(\frac{1}{\epsilon} \mathcal{R}_{(k)} + \partial_t \varphi_{(k)}^1 + \varphi_{(k)}^1 \cdot \nabla (\mathcal{U} + v^0 + \epsilon v^1) + \varphi_{(k)}^0 \cdot \nabla v^1\right) \nonumber \\
&\left.(v^1 \cdot \nabla \varphi_{(k)}^0 + v^0 \cdot \nabla \varphi_{(k)}^1) + (\varphi_{(k)}^1 \cdot \nabla \varphi^0 + \varphi_{(k)}^0 \cdot \nabla \varphi^1) + \epsilon (v^1 \cdot \nabla \varphi_{(k)}^1 + \varphi_{(k)}^1 \cdot \nabla \varphi^1)\right\}, && \nonumber
\end{align}
so that $\mathcal{L}^1 = \sum_{k=1}^N \mathcal{L}_{(k)}^1$.

We rewrite each term in $\mathcal{L}_{(k)}^1$ as follows:
\begin{align}
&\left(- \epsilon \Delta \varphi_{(k)}^1 + \mathcal{U} \cdot \nabla \varphi_{(k)}^1 - \frac{1}{\epsilon} \mathcal{M}_{(k)} + \nabla q_{(k)}^1\right)^i \nonumber\\
= &\left(- \epsilon \frac{\partial^2 \varphi_{(k)}^{1,i}}{\partial \xi_3^2} - \alpha^3 \frac{\partial \varphi_{(k)}^{1,i}}{\partial \xi_3} - \frac{1}{\epsilon} \mathcal{M}_{(k)}^i\right) \nonumber \\
&+ \left(- \epsilon S^i(\varphi_{(k)}^1) + \frac{\partial \mathcal{U}^3}{\partial \xi_3} \Big|_{\xi_3 = \zeta} \xi_3 \frac{\partial \varphi_{(k)}^{1,i}}{\partial \xi_3} + P^i(\mathcal{U}, \varphi_{(k)}^1) + (\nabla q_{(k)}^1)^i\right) - \epsilon T^i(\varphi_{(k)}^1) \nonumber\\
= &p_0(\epsilon,\xi_3) e^{-\frac{\alpha^3 \xi_3}{\epsilon}} + p_0(\epsilon,\xi_3) e^{-\frac{2 \alpha^3 \xi_3}{\epsilon}} + \sigma \epsilon, \qquad i=1,2, \nonumber
\end{align}
\begin{align}
&\left(- \epsilon \Delta \varphi_{(k)}^1 + \mathcal{U} \cdot \nabla \varphi_{(k)}^1 - \frac{1}{\epsilon} \mathcal{M}_{(k)} + \nabla q_{(k)}^1\right)^3 \nonumber \\
= &\left(- \epsilon \frac{\partial^2 \varphi_{(k)}^{1,3}}{\partial \xi_3^2} - \alpha^3 \frac{\partial \varphi_{(k)}^{1,3}}{\partial \xi_3} - \epsilon S^3(\varphi_{(k)}^1) + \frac{\partial \mathcal{U}^3}{\partial \xi_3} \Big|_{\xi_3 = \zeta} \xi_3 \frac{\partial \varphi_{(k)}^{1,3}}{\partial \xi_3} + P^3(\mathcal{U}, \varphi_{(k)}^1)\right) \nonumber\\
&+ \left((\nabla q^1)^3 - \frac{1}{\epsilon} \mathcal{M}_{(k)}^3\right) - \epsilon T^3(\varphi_{(k)}^1)\\
= &p_0(\epsilon,\xi_3) e^{-\frac{\alpha^3 \xi_3}{\epsilon}} + p_0(\epsilon,\xi_3) e^{-\frac{2 \alpha^3 \xi_3}{\epsilon}} + \sigma \epsilon, \nonumber
\end{align}
 where we have used \eqref{directional derivative}, \eqref{laplacian}, \eqref{excise higher order BL}, \eqref{higher order BL}, \eqref{excise higher order BL 2}, \eqref{grad pressure}), and 
\begin{align}
&\frac{1}{\epsilon} \mathcal{R}_{(k)} + \partial_t \varphi_{(k)}^1 + \varphi_{(k)}^1 \cdot \nabla (\mathcal{U} + v^0 + \epsilon v^1) + \varphi_{(k)}^0 \cdot \nabla v^1 \nonumber \\
= &\sum_{i=1,2,3} \left(p_0(\epsilon,\xi_3) e^{-\frac{\alpha^3 \xi_3}{\epsilon}} + p_0(\epsilon,\xi_3) e^{-\frac{2 \alpha^3 \xi_3}{\epsilon}} + \sigma \epsilon\right) e_{(k),i}, \nonumber
\end{align}
by \eqref{directional derivative}, \eqref{BL formula 2}, \eqref{BL formula 3}, \eqref{R_k}, \eqref{higher order BL}.

We note that both $v^0$ and $v^1$ satisfy $\gamma_- v^m \cdot n = 0$, and that $\varphi_{(k)}^0$ and $\varphi_{(k)}^1$ have analogous structures (with $\varphi_{(k)}^1$ containing additional $e^{-2\alpha^3 \xi_3/\epsilon}$ terms). Hence, proceeding as in \eqref{L^0 1} and \eqref{L^0 2} yields, we obtain  for $m,n=0,1$ and any $k,l$,
\begin{equation*}
v^m \cdot \nabla \varphi_{(k)}^n = \sum_{i=1,2,3} \left(p_0(\epsilon,\xi_3) e^{-\frac{\alpha^3 \xi_3}{\epsilon}} + p_0(\epsilon,\xi_3) e^{-\frac{2 \alpha^3 \xi_3}{\epsilon}} + \sigma \epsilon\right) e_{(k),i};
\end{equation*}
\begin{align}
&\varphi_{(k)}^m \cdot \nabla \varphi_{(l)}^n \nonumber \\
= &\sum_{i=1,2,3} \left(p_0(\epsilon,\xi_3) e^{-\frac{\alpha^3 \xi_3}{\epsilon}} + p_0(\epsilon,\xi_3) e^{-\frac{2 \alpha^3 \xi_3}{\epsilon}} + p_0(\epsilon,\xi_3) e^{-\frac{3 \alpha^3 \xi_3}{\epsilon}} + p_0(\epsilon,\xi_3) e^{-\frac{4 \alpha^3 \xi_3}{\epsilon}} + \sigma \epsilon^2\right) e_{(k),i}. && \nonumber
\end{align}
Combining the above formulas, we conclude that
\begin{equation*}
\mathcal{L}^1_{(k)} = \sum_{i=1,2,3} \sum_{r=1,2,3,4} \left(p_1(\epsilon,\xi_3) e^{-\frac{r \alpha^3 \xi_3}{\epsilon}} + \sigma \epsilon^2\right) e_{(k),i}.
\end{equation*}

Finally, we perform an energy estimate on \eqref{higher order remainder}. 
Proceeding similarly to what done in Section \ref{s:Sobolev}, we estimate
\begin{equation*}
\Big|\int \epsilon^2 (- \Delta v^1 + v^1 \cdot \nabla v^1) \cdot w\Big| \leq C \epsilon^4 + C \|w\|_2^2.
\end{equation*}
For $\mathcal{L}_{(k)}^1$, we have instead:
\begin{align}
\Big|\int \mathcal{L}_{(k)}^1 \cdot w\Big| = &\Big|\int \mathcal{L}_{(k)}^{1,i} w_i\Big| \leq \sum_{i=1,2,3} \sum_{r=1,2,3,4} \Big|\int (p_1(\epsilon,\xi_3) e^{-\frac{r \alpha^3 \xi_3}{\epsilon}} + \sigma \epsilon^2) w_i\Big| \nonumber \\
& \leq C \sum_{r=1,2,3,4} \left\|p_2(\epsilon,\xi_3) e^{-\frac{r \alpha_0 \xi_3}{\epsilon}}\right\|_2 \left\|\frac{w}{\xi_3}\right\|_2 + C \epsilon^2 \|w\|_2\\
&\leq \epsilon^{5/2} \|\nabla w\|_2 + C \epsilon^2 \|w\|_2 \nonumber\\
&\leq C \epsilon^4 + \frac{\epsilon}{4N} \|\nabla w\|_2^2 + C \|w\|_2^2.  \nonumber
\end{align}

For the nonlinear terms, we exploit estimates \eqref{nonlinear 1} and \eqref{nonlinear 2} so that
\begin{equation*}
\Big|\int (w \cdot \nabla (\mathcal{U} + v^0)) \cdot w\Big| \leq C \|w\|_2^2,
\end{equation*}
\begin{equation*}
\Big|\int (w \cdot \nabla \varphi^0_{(k)}) \cdot w\Big| \leq \frac{\epsilon}{4N} \|\nabla w\|_2^2 + C \|w\|_2^2,
\end{equation*}
and finally,
\begin{equation*}
\Big|\int (w \cdot \nabla \epsilon v^1) \cdot w\Big| \leq C \epsilon \|w\|_2^2,
\end{equation*}
\begin{align}
\Big|\int (w \cdot \nabla \epsilon \varphi_{(k)}^1) \cdot w\Big| &\leq \|\epsilon \nabla \varphi_{(k)}^1\|_\infty \|w\|_2^2 \nonumber \\
&\leq C \left\|p_0(\epsilon,\xi_3) e^{-\frac{\alpha^3 \xi_3}{\epsilon}} + p_0(\epsilon,\xi_3) e^{-\frac{2 \alpha^3 \xi_3}{\epsilon}} + \sigma \epsilon^2\right\|_\infty \|w\|_2^2 \nonumber \\
&\leq C \|w\|_2^2. \nonumber
\end{align}

Summing over k, we then obtain the energy estimate:
\begin{equation*}
\frac{d}{dt} \|w\|_2^2 + \epsilon \|\nabla w\|_2^2 \leq C \epsilon^4 + C \|w\|_2^2,
\end{equation*}
which by Gronwall's inequality implies 
\begin{equation}
\begin{cases}
\|w\|_{L^\infty L^2} \leq C \epsilon^2,\\
\|w\|_{L^2 H^1} \leq C \epsilon^{3/2}.
\end{cases}
\end{equation}

% final remark
\begin{Remark}
By iterating the above construction, it is possible to derive boundary layer correctors of arbitrarily high order under suitable higher-order compatibility conditions, thereby achieving convergence rates of arbitrarily high precision in the vanishing viscosity limit. More precisely, if one carries out an asymptotic expansion up to order $n \in \mathbb{N}_0$, leading to the approximate solution of order $n$
\begin{equation*}
\varphi_{(n)} = \varphi^0 + \sum_{j=1}^n (v^j + \varphi^j),
\end{equation*}
then we have the following error estimates:
\begin{equation}
\begin{cases}
\|v^\epsilon - v^0 - \varphi_{(n)}\|_{L^\infty L^2} \leq C \epsilon^{n+1},\\
\|v^\epsilon - v^0 - \varphi_{(n)}\|_{L^2 H^1} \leq C \epsilon^{n+\frac{1}{2}},
\end{cases}
\end{equation}
for some constant $C > 0$ independent of $\epsilon$, provided the compatibility condition
\begin{equation*}
\gamma_-(\partial_t^m u^0)|_{t=0} \cdot \tau = (\partial_t^m \alpha)|_{t=0} \cdot \tau,
\end{equation*}
or equivalently,
\begin{equation*}
\gamma_-(\partial_t^m v^0)|_{t=0} \cdot \tau = 0,
\end{equation*}
is satisfied for $m = 0, \dots, n$.
\end{Remark}

\bigskip

\section*{Acknowledgments}
A. L. Mazzucato was supported in part by NSF grants DMS-2206453 and DMS-2511023.
D. Wang was supported in part by NSF grants DMS-2219384 and DMS-2510532.
W. Wei was supported in part by NSF grant DMS-2219384.

\bibliography{InflowOutflow}
\bibliographystyle{abbrv}

\end{document}